\input amstex
\input amsppt.sty
\magnification=\magstep1
\hsize=32truecc
\vsize=22.5truecm
\baselineskip=16truept
\NoBlackBoxes
\TagsOnRight \pageno=1 \nologo
\def\Z{\Bbb Z}
\def\N{\Bbb N}

\def\l{\left}
\def\r{\right}
\def\bg{\bigg}
\def\({\bg(}
\def\[{\bg\lfloor}
\def\){\bg)}
\def\]{\bg\rfloor}
\def\t{\text}
\def\f{\frac}

\def\bi{\binom}
\def\eq{\equiv}

\def\ls{\leqslant}
\def\gs{\geqslant}
\def\mo{\roman{mod}}

\def\ve{\varepsilon}
\def\al{\alpha}
\def\da{\delta}

\def\Proof{\noindent{\it Proof}}

\def\Remark{\medskip\noindent{\it  Remark}}

\def\Ack{\medskip\noindent {\bf Acknowledgments}}
\hbox {J. Number Theory 133(2013), no.\,8, 2794--2812.}
\bigskip
\topmatter
\title On functions taking only prime values\endtitle
\author Zhi-Wei Sun\endauthor
\leftheadtext{Zhi-Wei Sun}
 \rightheadtext{On functions taking only prime values}
\affil Department of Mathematics, Nanjing University\\
 Nanjing 210093, People's Republic of China
  \\  zwsun\@nju.edu.cn
  \\ {\tt http://math.nju.edu.cn/$\sim$zwsun}
\endaffil
\abstract For $n=1,2,3,\ldots$ define $S(n)$ as the smallest integer
$m>1$ such that those $2k(k-1)$ mod $m$ for $k=1,\ldots,n$ are
pairwise distinct; we show that $S(n)$ is the least prime greater
than $2n-2$ and hence the value set of the function $S(n)$ is
exactly the set of all prime numbers. For every $n=4,5,\ldots$, we
prove that the least prime $p>3n$ with $p\eq1\ (\mo\ 3)$ is just the
least positive integer $m$ such that $18k(3k-1)\ (k=1,\ldots,n)$ are
pairwise distinct modulo $m$. For $d\in\{4,6,12\}$ and
$n=3,4,\ldots$, we show that the least prime $p\gs 2n-1$ with
$p\eq-1\ (\mo\ d)$ is the smallest integer $m$ such that those
$(2k-1)^d$ for $k=1,\ldots,n$ are pairwise distinct modulo $m$.
We also pose several challenging conjectures on primes. For example, we
find a surprising recurrence for primes, namely, for every
$n=10,11,\ldots$ the $(n+1)$-th prime $p_{n+1}$ is just the least
positive integer $m$ such that $2s_k^2\ (k=1,\ldots,n)$ are pairwise
distinct modulo $m$ where $s_k=\sum_{j=1}^k(-1)^{k-j}p_j$. We also
conjecture that for any positive integer $m$ there are consecutive
primes $p_k,\ldots,p_n\ (k<n)$ not exceeding $2m+2.2\sqrt m$ such
that $m=p_n-p_{n-1}+\cdots+(-1)^{n-k}p_k$.
\endabstract
\thanks 2010 {\it Mathematics Subject Classification}.\,Primary 11A41, 11Y11;
Secondary 05A10, 11A07, 11B75, 11N13.
\newline\indent Supported by the National Natural Science
Foundation (grant 11171140) of China.
\endthanks
\endtopmatter
\document

\heading{1. Introduction}\endheading

Prime numbers play a central role in number theory (see the excellent book [CP] on primes by R. Crandall and C. Pomerance). It is
known that there is no non-constant polynomial with integer
coefficients, even in several variables, which takes only prime
values. Many mathematicians ever tried in vain to find a nontrivial
number -theoretic function whose values are always primes. In 1947 W. H. Mills [Mi]
showed that there exists a real number $A$ such that $\lfloor
A^{3^n}\rfloor$ is prime for any $n=1,2,3,\ldots$;
unfortunately such a constant $A$ cannot be effectively found.

In 2012, the author conjectured that for any $m=12,13,\ldots$ the largest positive integer $n$
with $\bi{2k}k\ (k=1,\ldots,n)$ pairwise distinct modulo $m$ does not exceed $0.6\sqrt m\log m$.
Motivated by this, he made the following conjecture.

\proclaim{Conjecture 1.1} {\rm (i) ([S12a])} For $n\in\Z^+=\{1,2,3,\ldots\}$ define $s(n)$
 as the smallest integer $m>1$ such that
$$\bi{2k}k\ \ (k=1,\ldots,n)$$
are pairwise distinct modulo $m$. Then
all those $s(1),s(2),\ldots$ are primes!

{\rm (ii) ([S12b])} For $n\in\Z^+$ let $t(n)$ denote the least integer $m>1$ such that
$$|\{k!\ \mo\ m:\ k=1,\ldots,n\}|=n.$$
Then $t(n)$ is prime with the only exception $t(5)=10$.
\endproclaim

The author verified both parts of Conjecture 1.1 for $n\ls 2000$. Later, Laurent Bartholdi and Qing-Hu Hou verified
parts (i) and (ii) of Conjecture 1.1 for all $n\in[2001,5000]$ and $n\in[2001,10000]$ respectively.

In 1985 L. K. Arnold, S. J. Benkoski and B. J. McCabe [ABM] defined $D(n)$ for $n\in\Z^+$ as
the smallest positive integer $m$ such that $1^2,2^2,\ldots,n^2$ are pairwise distinct modulo $m$, and they showed that
if $n>4$ then $D(n)$ is the smallest integer $m\gs 2n$ such that $m$ is $p$ or $2p$ with $p$ an odd prime.
(Note that if $p=2q+1$ with $p$ and $q$ both prime then the prime $p$ is not contained in the range $\{D(n):\ n\in\Z^+\}$.)
This stimulated later studies of
characterizing
$$D_f(n):=\min\{m\in\Z^+:\ f(1),f(2),\ldots,f(n)\ \t{are distinct modulo}\ m\}$$
for some special polynomials $f(x)\in\Z[x]$ including powers of $x$ and Dickson polynomials of degrees relatively prime to $6$
(see, e.g., [BSW, MM, Z] and the references therein).
However, the value sets of those $D_f$ considered in papers along this line are usually somewhat complicated
and they contain infinitely many composite numbers. Note also that $D_f(1)$ is just $1$, not a prime.

Now we present a simple function whose set of values is exactly the set of all prime numbers.

\proclaim{Theorem 1.1} {\rm (i)} For $n\in\Z^+$ let $S(n)$ denote the smallest
integer $m>1$ such that those $2k(k-1)$ mod $m$ for
$k=1,\ldots,n$ are pairwise distinct. Then $S(n)$ is the least prime greater than $2n-2$.

{\rm (ii)} For $n\in\Z^+$ let $T(n)$ denote the least
integer $m>1$ such that those $k(k-1)$ mod $m$ with
$1\ls k\ls n$ are pairwise distinct. Then we have
$$T(n)=\min\{m\gs 2n-1:\ m\ \t{is a prime or a positive power of}\ 2\}.\tag1.1$$
\endproclaim
\Remark\ 1.1. (a) The way to generate all primes via Theorem 1.1(i) is simple in concept, but it has no advantage in algorithm.
Nevertheless, Theorem 1.1(i) is of certain theoretical interest since it provides a surprising new characterization of primes.

(b) By modifying our proof of Theorem 1.1(i), we are also able to show that for any $d,n\in\Z^+$
with $n\gs \lfloor d/2\rfloor+4$ the least prime $p\gs 2n+d$
is just the smallest $m\in\Z^+$ such that $2k(k+d)\ (k=1,\ldots,n)$ are pairwise distinct modulo $m$.
(Similar results hold for $d\in\{0,-2\}$ and $n\in\{5,6,\ldots\}$.)
\medskip

Below are four more related theorems.

\proclaim{Theorem 1.2} {\rm (i)} For any positive integer $n$, the number $2^{\lceil\log_2n\rceil}$ (the least power of two not smaller than $n$)
is the least positive integer $m$ such that
those $k(k-1)/2\ (k=1,\ldots,n)$ are pairwise distinct modulo $m$.

{\rm (ii)} Let $d\in\{2,3\}$ and $n\in\Z^+$. Take
the smallest positive integer $m$ such that $|\{k(dk-1)\ \mo\ m:\ k=1,\ldots,n\}|=n$.
Then $m$ is the least power of $d$ not smaller than $n$, i.e., $m=d^{\lceil\log_d n\rceil}$.

{\rm (iii)} Let $n\in\{4,5,\ldots\}$ and take the least positive integer $m$ such that $18k(3k-1)\ (k=1,\ldots,n)$ are pairwise distinct modulo $m$.
Then $m$ is the least prime $p>3n$ with $p\eq1\ (\mo\ 3)$.
\endproclaim
\Remark\ 1.2. We are also able to prove some other results similar to those in Theorem 1.2. For example,
for each $n=5,6,7,\ldots$ the first prime $p\eq-1\ (\mo\ 3)$ after $3n$
is just the least $m\in\Z^+$ such that those $18k(3k+1)\ (k=1,\ldots,n)$
are pairwise distinct modulo $m$. Also, if $f(n)$ denotes the least $m\in\Z^+$ with $|\{4k(4k-1)\ \mo\ m:\ k=1,\ldots,n\}|=n$
and $g(n)$ denotes the least $m\in\Z^+$ with $|\{4k(4k+1)\ \mo\ m:\ k=1,\ldots,n\}|=n$, then
$f(n)$ with $n\gs5$ is the least prime $p>(8n-4)/3$ with $p\eq1\ (\mo\ 4)$,
and $g(n)$ with $n\gs6$ is the least prime $p>(8n-2)/3$ with $p\eq-1\ (\mo\ 4)$.

\proclaim{Theorem 1.3} For $d,n\in\Z^+$ let $\lambda_d(n)$ be the smallest integer $m>1$ such that those
$(2k-1)^d\ (k=1,\ldots,n)$ are pairwise distinct modulo $m$. Then $\lambda_d(n)$ with $d\in\{4,6,12\}$ and $n>2$
is the least prime $p\gs 2n-1$ with $p\eq-1\ (\mo\ d)$.
\endproclaim

\proclaim{Theorem 1.4} Let $q$ be an odd prime. Then the smallest integer $m>1$ such that those $k^q(k-1)^q$
$(k=1,\ldots,n)$ are pairwise distinct modulo $m$, is just the least prime $p\gs 2n-1$ with $p\not\eq1\ (\mo\ q).$
\endproclaim

\proclaim{Theorem 1.5} Define $s_n=\sum_{k=1}^n(-1)^{n-k}p_k$ for all $n\in\Z^+$, where $p_k$ denotes the $k$-th prime.
Then, for any $n\in\Z^+$ those $2s_k^2\ (k=1,\ldots,n)$ are pairwise distinct modulo $p_{n+1}$.
\endproclaim
\Remark\ 1.3. All terms of the sequence $s_1,s_2,s_2,\ldots$ are positive integers.
In fact, if $n\in\Z^+$ is even then $s_n=\sum_{k=1}^{n/2}(p_{2k}-p_{2k-1})>0$; if $n\in\Z^+$ is odd then
$s_n=\sum_{k=1}^{(n-1)/2}(p_{2k+1}-p_{2k})+p_1>0$.
Here we list the values of $s_1,\ldots,s_{15}$.
$$\gather s_1=2,\ s_2=1,\ s_3=4,\ s_4=3,\ s_5=8,\ s_6=5,\ s_7=12,\ s_8=7,
\\ s_9=16,\ s_{10}=13,\ s_{11}=18,\ s_{12}=19,\ s_{13}=22,\ s_{14}=21,\ s_{15}=26.
\endgather$$
The sequence $0,s_1,s_2,\ldots$ was first introduced by N.J.A. Sloane and J. H. Conway [SC].
We conjecture that for any integers $m>0$ and $r$ there are infinitely many $n\in\Z^+$ with $s_n\eq r\ (\mo\ m)$.

\medskip

In the next section we will present two auxiliary theorems.
Section 3 is devoted to our proofs of Theorems 1.1 and 1.2.
In Section 4 we will show Theorems 1.3-1.5.

Motivated by Theorem 1.5 we raise the following conjecture on recurrence for primes
which allows us to compute $p_{n+1}$ in terms of $p_1,\ldots,p_n$.

\proclaim{Conjecture 1.2} Let $n\in\Z^+$ with $n\not=1,2,4,9$. Then $p_{n+1}$ is the smallest positive integer $m$
such that those $2s_k^2\ (k=1,\ldots,n)$ are pairwise distinct modulo $m$.
\endproclaim
\Remark\ 1.4. (a) We have verified Conjecture 1.2 for all $n\ls 10^5$. Note that 9 is the least $m\in\Z^+$
with $2s_1^2,2s_2^2,2s_3^2,2s_4^2$ pairwise distinct modulo $m$, and 25 is the least $m\in\Z^+$ with
$|\{2s_k^2\ \mo\ m: \ k=1,\ldots,9\}|=9$.

(b) Define $b(n)$ as the least power of two
modulo which $s_1,\ldots,s_n$ are pairwise incongruent.
We conjecture that $b(n)$ is the least $m\in\Z^+$ such that
$2s_k^2-s_k\ (k=1,\ldots,n)$ are pairwise distinct modulo $m$,
and moreover $\{b(n):n\in\Z^+\}=\{2^a:\ a=0,1,2,\ldots\}$.
\medskip

Inspired by Conjecture 1.2, we find the following surprising conjecture on representations of integers by alternating sums of consecutive primes.

\proclaim{Conjecture 1.3} For any positive integer $m$, there are consecutive primes $p_k,\ldots,p_n\ (k<n)$ not exceeding $2m+2.2\sqrt m$
such that
$$m=p_n-p_{n-1}+\cdots+(-1)^{n-k}p_k.$$
\endproclaim
\Remark\ 1.5. We also conjecture that $2m+2.2\sqrt m$ in Conjecture 1.3 can be replaced by $m+4.6\sqrt m$ if $m$ is odd.
If the upper bound $2m+2.2\sqrt{m}$ is replaced by $3m$, then we may require additionally that
$p_k-1$ and $p_n+1$ are both practical numbers (cf. [S13]).
We have verified Conjecture 1.3 for $m=1,\ldots,10^5$. To
illustrate the conjecture, we look at a few concrete examples:
$$\gather1=3-2,\ \ 2=5-3,\ \ 3=7-5+3-2,\ \ 4=11-7,\ \ 5=7-5+3,
\\ 8=11-7+5-3+2,\ \ 11=19-17+13-11+7,
\\ 20=41-37+31-29+23-19+17-13+11-7+5-3,
\\ 303=p_{76}-p_{75}+\cdots-p_{53}+p_{52}\ \ \t{with}\ p_{76}=383=303+\lfloor 4.6\sqrt{303}\rfloor,
\\ 2382=p_{652}-p_{651}+\cdots+p_{44}-p_{43}\ \ \t{with}\ p_{652}=4871= 2\cdot 2382+\lfloor 2.2\sqrt{2382}\rfloor.
\endgather$$
The author would like to offer 1000 US dollars as the prize for the first correct proof of Conjecture 1.3.
We also have some other conjectures on representations involving alternating sums of consecutive primes, for example, every $m=3,4,\ldots$
can be written in the form $p+s_n$, where $p$ is a Sophie Germain prime and $n$ is a positive integer.
\medskip

We also have a conjecture involving sums of consecutive primes.

\proclaim{Conjecture 1.4}  For $k\in\Z^+$ let $S_k$ denote the sum of the first $k$ primes $p_1,\ldots,p_k$.

{\rm (i)} For $n\in\Z^+$ define $S^+(n)$ as the least integer $m>1$ such that $m$ divides none of $S_i!+S_j!$ with $1\ls i<j\ls n$.
Then $S^+(n)$ is always a prime, and $S^+(n)<S_n$ for every $n=2,3,4,\ldots$.

{\rm (ii)} For $n\in\Z^+$ define $S^-(n)$ as the least integer $m>1$ such that $m$ divides none of those $S_i!-S_j!$ with $1\ls i<j\ls n$.
Then $S^-(n)$ is always a prime, and $S^-(n)<S_n$ for every $n=2,3,4,\ldots$.

{\rm (iii)} For any positive integer $n$ not dividing $6$, the least integer $m>1$ such that $2S_k^2\ (k=1,\ldots,n)$ are pairwise distinct modulo $m$
is a prime smaller than $n^2$.
\endproclaim
\Remark\ 1.6. When $n>1$, clearly $S_n!\pm S_{n-1}!\eq0\ (\mo\ m)$ for any $m=1,\ldots,S_{n-1}$, and hence both $S^+(n)$ and $S^-(n)$
are greater than $S_{n-1}$. Thus, by the conjecture we should have $S^+(n)<S_n<S^+(n+1)$ and $S^-(n)<S_n<S^-(n+1)$  for all $n=2,3,\ldots$.
Conjecture 1.4 implies that for any $n=2,3,\ldots$ the interval $(S_{n-1},S_{n})$ contains the primes $S^+(n)$ and $S^-(n)$,
which are actually very close to $S_{n-1}$.
However, it seems very challenging to prove that $(S_n,S_{n+1})$ contains a prime for any $n\in\Z^+$. Note that
$$S_n\sim\sum_{k=1}^n k\log k\sim\int_1^n x\log xdx=\f{x^2}2\log x\bigg|_1^n-\int_1^n\f{x^2}2(\log x)'dx\sim\f{n^2}2\log n$$
as $n\to+\infty$, and the Legendre conjecture asserts that the interval $(n^2,(n+1)^2)$ contains a prime for any $n\in\Z^+$.
We conjecture that the number of primes in the interval $(S_n,S_{n+1})$ is asymptotically equivalent to $cn/2$ as $n\to+\infty$,
where $c\gs1$ is a constant
(whose value is probably 1).
\medskip

Our following conjecture allows us to produce primes via products of consecutive primes.

\proclaim{Conjecture 1.5} For $k\in\Z^+$ let $P_k$ denote the product of the first $k$ primes $p_1,\ldots,p_k$.

{\rm (i)} For $n\in\Z^+$ define $w_1(n)$ as the least integer $m>1$ such that $m$ divides none of those $P_i-P_j$ with $1\ls i<j\ls n$.
Then $w_1(n)$ is always a prime.

{\rm (ii)} For $n\in\Z^+$ define $w_2(n)$ as the least integer $m>1$ such that $m$ divides none of those $P_i+P_j$ with $1\ls i<j\ls n$.
Then $w_2(n)$ is always a prime.

{\rm (iii)} We have $w_1(n)<n^2$ and $w_2(n)<n^2$ for all $n=2,3,4,\ldots$.
\endproclaim

\Remark\ 1.7. (a) Clearly $w_i(n)\ls w_i(n+1)$ for $i=1,2$ and $n\in\Z^+$.
Since $P_1,\ldots,P_n$ are pairwise distinct modulo $w_1(n)$, we have $w_1(n)\gs n$
and hence $W_1=\{w_1(n):\ n\in\Z^+\}$ is an infinite set. For any integer $m>1$, there is an odd prime $p_n\eq-1\ (\mo\ m)$ and hence
$P_{n-1}+P_n=P_{n-1}(1+p_n)\eq0\ (\mo\ m)$. Thus $W_2=\{w_2(n):\ n\in\Z^+\}$ is also infinite. If $w_i(n)=p_k$, then $k\gs n$ since
$P_k\pm P_{k+1}\eq0\ (\mo\ p_k)$. Thus it follows from Conjecture 1.5(ii) that $w_2(n)>n$ for all $n\in\Z^+$, in other words,
for each $n=2,3,4,\ldots$ there are $1\ls j<k\ls n$ such that $P_j+P_k\eq0\ (\mo\ n)$.
For $n=2,3,4,\ldots$ we conjecture further that $P_n\eq P_j\eq-P_k\ (\mo\ n)$ for some $j,k\in\{1,\ldots,n-1\}$.
This seems simple but we are unable to prove it.

(b) The author [S12c] listed values of $w_1(n)$ for $n=1,\ldots,1172$, and values of $w_2(n)$ for $n=1,\ldots,258$. Later W. B. Hart [H]
reported that he had verified Conjecture 1.5 for all $n\ls 10^5$.
\medskip

A prime is said to be of the first kind (or the second kind) if it belongs to $W_1=\{w_1(n):\ n\in\Z^+\}$ (or $W_2=\{w_2(n):\ n\in\Z^+\}$, resp.).
Here we list the first 20 primes of each kind.

Primes of the first kind: 2, 3, 5, 11, 23, 29, 37, 41, 47, 73, 131,
151, 199, 223, 271, 281, 353, 457, 641, 643, $\ldots$

Primes of the second kind: 2, 3, 5, 7, 11, 19, 23, 47, 59, 61, 71, 101,
113, 223, 487, 661, 719, 811, 947, 1327, $\ldots$
\medskip

The famous Artin conjecture for primitive roots states that if an integer $a$ is neither $-1$ nor
a square then there are infinitely many primes $p$ having $a$ as a primitive root modulo $p$.
This is open for any particular value of $a$.
Concerning Artin's conjecture the reader may consult the excellent survey of R. Murty [Mu] and the book [IR, p.\,47].
In Section 5 we will present more conjectures which are similar to Conjecture 1.1 or related to the Artin conjecture.

\heading{2. Two auxiliary theorems}\endheading

\proclaim{Theorem 2.1} Let $m>1$ and $n>1$ be integers such
that those $k(k-1)$ for $k=1,\ldots,n$ are pairwise distinct
modulo $m$.

{\rm (i)} We have $m\gs 2n-1$.

{\rm (ii)} If $n\gs15$ and $m\ls 2.4n$, then $m$ is a prime or a power of two.
\endproclaim

\medskip\noindent{\it Proof of Theorem} {\rm 2.1(i)}. Suppose on the contrary that $m\ls 2n-2$.
Then $n\gs m/2+1$. If $m$ is even, then
$$\l(\f m2+1\r)\l(\f m2+1-1\r)-\f m2\l(\f m2-1\r)=m\eq0\pmod m.$$
If $m$ is odd, then $(m+3)/2\ls n$ and
$$\f{m+3}2\l(\f{m+3}2-1\r)-\f{m-1}2\l(\f{m-1}2-1\r)=2m\eq0\pmod m.$$
So we get a contradiction as desired. \qed

The next task in this section is to prove Theorem 2.1(ii).
In the following two lemmas, we fix $n\gs15$ and $m\in[2n-1,2.4n]$ and assume that those $k(k-1)$ mod $m$ $(1\ls k\ls n)$ are pairwise distinct.

\proclaim{Lemma 2.1} $m\not= 2p$ for any odd
prime $p$.
\endproclaim
\Proof. Suppose that $m=2p$ with $p$ an odd prime. Note that
$$\f{p+3}2\l(\f{p+3}2-1\r)-\f{p-1}2\l(\f{p-1}2-1\r)=2p\eq0\pmod{2p}.$$
and hence $(p+3)/2>n$. So $2n-1\ls p=m/2\ls 1.2n$, which is impossible. \qed

\proclaim{Lemma 2.2} $p^2\nmid m$ for any odd prime $p$.
\endproclaim
\Proof. Suppose that $m=p^2q$ with $p$ an odd prime and $q\in\Z^+$.
Set $k=(p+1)/2$ and $l=k+pq\ls 2pq$. Then
$$l(l-1)-k(k-1)=(l-k)(l+k-1)=pq(pq+2k-1)\eq0\pmod{p^2q}$$
and hence we must have $2pq>n$. If $p>3$, then
$$n<\f{2m}p\ls\f 25m\ls \f25\times 2.4n<n$$
which is impossible. When $p=3$, we also have a contradiction since $l=2+3q=2+m/3\ls 2+0.8n\ls n$.
\qed

\medskip
\noindent {\it Proof of Theorem} {\rm 2.1(ii)}. Suppose that $n\gs15$ and $m\ls 2.4n$.
We want to deduce a contradiction under the assumption that
$m$ is neither a prime nor a power of two.

By Lemmas 2.1 and 2.2, we may write $m=pq$ with $p$ an odd prime, $q>2$ and $p\nmid q$.

Take an integer $k\in[1,q/(2,q)]$ such that
$$k\eq\f{1-p}2\ \l(\mo\ \f q{(2,q)}\r),$$
where $(2,q)$ is the greatest common divisor of $2$ and $q$.
Set $l=k+p$. Then
$$l(l-1)-k(k-1)=p(2k-1+p)\eq0\ \ (\mo\ pq).$$
If $2\mid q$, then $q\gs 4$ and hence
$$l\ls p+\f q2=\f{m}{q}+\f{m}{2p}\ls\f{m}4+\f{m}6=\f 5{12}m\ls\f 5{12}\times 2.4n=n$$
which contradicts the property of $m$. Thus $2\nmid q$ and
$$l\ls p+q=\f{m}p+\f{m}q\ls\l(\f 1p+\f1q\r)2.4n.$$
 If both $p$ and $q$ are greater than 3,  then
 $$\f1{p}+\f1{q}\ls\f25<\f 5{12}$$
 and hence $l<\f5{12}2.4n=n$ which leads to a contradiction.
 So $m$ cannot have two distinct prime divisors greater than 3.
In view of Lemma 2.2, we may assume that $m=pq$ with $q=3$.
Note that
 $$l\ls p+q=\f{m}3+3\ls \f{2.4n}3+3=0.8n+3\ls n$$
since $n\gs15$. So we get a contradiction.

 In view of the above, we have completed the proof of Theorem 2.1. \qed

\proclaim{Theorem 2.2} Let $n>1$ and $m\gs2n-1$ be integers.

{\rm (i)} Suppose that $m$ is a prime or a power of two. Then
 $k(k-1)\not\eq l(l-1)\pmod m$ for any $1\ls k<l\ls n$.

{\rm (ii)} If $m$ is a power of two not exceeding $2.4n$, then
 $2k(k-1)\eq2l(l-1)\pmod m$ for some $1\ls k<l\ls n$.
\endproclaim
\Proof. (i) To prove part (i) we distinguish two cases.

{\it Case} 1. $m=2^a$ for some $a\in\Z^+$.

In this case, $n\ls(m+1)/2=2^{a-1}+1/2$ and hence $n\ls 2^{a-1}$. For
any $1\ls k<l\ls n$, we have $0<l-k<n\ls 2^{a-1}$ and $0<l+k-1<2n\ls
2^a$, hence $$l(l-1)-k(k-1)=(l-k)(l+k-1)\not\eq0\pmod{2^a}$$ since
one of $l-k$ and $l+k-1$ is odd.

{\it Case} 2. $m$ equals an odd prime $p$.

If $1\ls k<l\ls n$, then $0<l-k<n\ls (p+1)/2<p$ and $l+k-1<2n-1\ls p$,
therefore
$$l(l-1)-k(k-1)=(l-k)(l+k-1)\not\eq0\pmod p.$$

(ii) As $2k(k-1)\eq0\ (\mo\ 4)$ for any $k=1,\ldots,n$, we just assume that $m=2^a$ with $a>2$.
Take $k=2^{a-2}$ and $l=k+1$. Then
$$2l(l-1)-2k(k-1)=2(2^{a-2}+1)2^{a-2}-2\times 2^{a-2}(2^{a-2}-1)=2^a\eq0\pmod{2^a}$$
and $k<l=2^{a-2}+1<2^a/2.4\ls n$.

Combining the above we have completed the proof. \qed

\heading{3. Proofs of Theorems 1.1 and 1.2}\endheading

As in Theorem 1.1, let $S(n)$ (or $T(n)$) denote the least integer $m>1$ such that those $2k(k-1)$ (or $k(k-1)$, resp.)
for $k=1,\ldots,n$ are pairwise distinct modulo $m$.

\proclaim{Lemma 3.1} For any positive integer $n$ we have $2n-1\ls T(n)\ls S(n)\ls 2.4n$.
\endproclaim
\Proof. The case $n=1$ is trivial since $S(1)=T(1)=2$. Below we assume $n\gs2$.

As those $2k(k-1)\ (k=1,\ldots,n)$ are pairwise distinct modulo $S(n)$,
those $k(k-1)\ (k=1,\ldots,n)$ are also pairwise distinct modulo $S(n)$ and hence $S(n)\gs T(n)$.
Note that $T(n)\gs 2n-1$ by Theorem 2.1(i).

By J. Nagura [N], for $m=25,26,\ldots$ the interval $[m, 1.2m]$ contains a prime.
Thus, if $n\gs 13$ then there is a prime in the interval $[2n-1,2.4n]$. For $n=2,\ldots,12$
we can easily check that the interval $[2n-1,2.4n]$ does contain primes.
By P. Dusart [D, Section 4], for $x\gs 3275$ there is a prime $p$ such that
$$x\ls p\ls x\l(1+\f1{2\log^2x}\r)\ls x\l(1+\f1{2\log^23275}\r)<1.01x;$$
this provides another way to show that $[2n-1,2.4n]$ contains at least a prime.
So there exists an odd prime $p\in[2n-1,2.4n]$ and hence $S(n)\ls p\ls 2.4n$ by Theorem 2.2(i).
(For $1\ls k<l\ls n$, clearly $k(k-1)\not\eq l(l-1)\ (\mo\ p)$ if and only if $2k(k-1)\not\eq2l(l-1)\ (\mo\ p)$.)
We are done. \qed

\medskip
\noindent{\it Proof of Theorem} 1.1. We want to prove that $S(n)$ is the least prime greater than $2n-2$ and $T(n)$ is the least integer $m\gs2n-1$
with $m$ a prime or a positive power of $2$. For $n=1,\ldots,14$ these can be easily verified.

Now assume that $n\gs15$. By Lemma 3.1, Theorem 2.1(ii) and Theorem 2.2(ii), $S(n)$ must be an odd prime in the interval $[2n-1,2.4n]$.
In view of Theorem 2.2(i), $S(n)$ is the least prime greater than $2n-2$.

By Lemma 3.1, $T(n)\in[2n-1,2.4n]$.
Applying Theorem 2.1(ii) we see that $T(n)$ is either a prime or a power of two.
Combining this with Theorem 2.2(i) we immediately get (1.1). \qed

\medskip
\noindent{\it Proof of Theorem} {\rm 1.2(i)}. Let $n\in\Z^+$ and
take the smallest positive integer $m$ such that those $k(k-1)/2\ (1\ls k\ls n)$ are pairwise distinct modulo $m$.
We want to prove that $m=2^h$ where $h:=\lceil\log_2n\rceil$. This is trivial when $n=1$.

Below we let $n>1$ and hence $h>0$. Note that $2^{h-1}<n\ls 2^h$.

Clearly $m\gs n$. As $2^{h+1}>2n-1$, by Theorem 2.2(i), those $k(k-1)$ $(k=1,\ldots,n)$
are pairwise distinct modulo $2^{h+1}$. It follows that $m\ls 2^h<2n$. If $m$ is odd, then $m\ls 2n-3$ and
$$\f12\cdot\f{m+3}2\l(\f{m+3}2-1\r)-\f12\cdot\f{m-1}2\l(\f{m-1}2-1\r)=m\eq0\pmod m.$$
So $m$ must be even.

Suppose that $m\not=2^h$. Then $m$ has the form $2p^aq$ with $p$ an odd prime, $a,q\in\Z^+$ and $p\nmid q$.
Let $k$ be the least positive residue of $(1-p^a)/2$ mod $2q$ and set $l=k+p^a$. Observe that
$$l(l-1)-k(k-1)=(l-k)(l+k-1)=p^a(2k-1+p^a)\eq0\ (\mo\ 4p^aq)$$
and thus $l(l-1)/2\eq k(k-1)/2\ (\mo\ m)$. Clearly,
$$l\ls 2q+p^a=\f m{p^a}+\f m{2q}<\l(\f{2}{p^a}+\f1q\r)n.$$
Thus we must have
$$\f2{p^a}+\f1q>1$$ and hence $q<3$. Thus $m=2p^a$ or $m=4\times3=12$.
When $n\ls12$ we can easily check that $m\not=12$. For $n>12$ we have $m\gs n>12$.
Therefore $m=2p^a$.

Note that $m/2+1\ls n$. If  $p^a\eq1\ (\mo\ 4)$, then
$$\f{p^a(p^a-1)}2-\f{1(1-1)}2=2p^a\f{p^a-1}4\eq0\ (\mo\ 2p^a)$$
and $p^a=m/2<n$;
if $p^a\eq3\ (\mo\ 4)$, then
$$\f{(p^a+1)p^a}2-\f{1(1-1)}2=2p^a\f{p^a+1}4\eq0\ (\mo\ 2p^a)$$
and $p^a+1=m/2+1\ls n$. So we get a contradiction.

The proof of Theorem 1.2(i) is now complete. \qed

\medskip
\noindent{\it Proof of Theorem} {\rm 1.2(ii)}. Fix $d\in\{2,3\}$ and $n\in\Z^+$, and take the least $m\in\Z^+$ such that
those $k(dk-1)\ (k=1,\ldots,n)$ are pairwise distinct modulo $m$. We want to prove that $m=d^{\lceil\log_d n\rceil}$.
This can be easily verified in the case $n\ls 7$.

Below we assume $n>7$ and hence $m\gs n\gs 8$. Suppose that $d^{h-1}<n\ls d^h$ where $h\in\Z^+$.
For $1\ls k<l\ls n$, clearly $0<l-k<n\ls d^h$ and hence
$$l(dl-1)-k(dk-1)=(l-k)(d(l+k)-1)\not\eq0\pmod{d^h}.$$
Thus $m\ls d^h<dn$.

When $m\eq -1\ (\mo\ d)$, we have $1<l=(m+1)/d-1<(m+1)/d\ls n$ and
$$l(dl-1)-1(d\cdot1-1)=\l(\f{m+1}d-1\r)((m+1-d)-1)-(d-1)\eq0\pmod m,$$
which contradicts the choice of $m$. So we have $m\not\eq -1\ (\mo\ d)$.
When $d=3$ and $m\eq1\ (\mo\ d)$, for $k=(m-1)/3$ and $l=(m+2)/3$, we have
$1\ls k<l\ls n$ and
$$l(dl-1)-k(dk-1)=(l-k)(d(l+k)-1)=2m\eq0\pmod m,$$
which also contradicts the choice of $m$. Therefore $m\not\eq\pm 1\ (\mo\ d)$
and hence $d\mid m$.

Write $m=d^aq$ with $a,q\in\Z^+$ and $d\nmid q$. Set $\da=d^a-\ve_q$, where
$$\ve_q=\cases -(\f{-1}q)&\t{if}\ d=2\ \t{and}\ a=1,
\\(\f{-1}q)&\t{if}\ d=2\ \t{and}\ a\gs2,
\\(\f{q}3)&\t{if}\ d=3,\endcases$$
and $(-)$ denotes the Legendre symbol.
Note that
$$\f{\da q+1}d=d^{a-1}q-\f{\ve_qq-1}d\in\Z\ \ \t{and}\ \ \f{\da q+1}d\eq d^a\ (\mo\ 2).$$
Thus both
$$k=\f12\l(\f{\da q+1}d-d^a\r)\ \ \t{and}\ \ l=\f12\l(\f{\da q+1}d+d^a\r)$$
are integers, and
$$l(dl-1)-k(dk-1)=(l-k)(d(l+k)-1)=d^a(\da q)\eq0\ (\mo\ m).$$
As
$$\f{\da q+1}d+d^a\ls d^{a-1}q+\f{q+1}d+d^a=\f{m+1}{d}+\f{m}{d^{a+1}}+\f mq$$
and $m<dn$, we have
$$2l<n+\f n{d^a}+\f{dn}q.$$

{\it Case} 1. $q\gs d+1$.

As $6(2\cdot6-1)-2(2\cdot2-1)=60=5(3\cdot5-1)-2(3\cdot2-1)$, we have $m\nmid 60$.
If $a=1$, then $q>5$, hence
$$\f{\da q+1}d\gs \f{(d-1)q+1}d>\f{5(d-1)-1}d=d$$
and
$$2l<n+\f nd+\f{dn}q\ls n+\f n2+\f{3n}6=2n.$$
When $a\gs2$, we have $q>d+(d-1)/(d^a-1)$ and hence
$$\f{\da q+1}d=d^{a-1}q-\f{\ve_qq-1}d\gs d^{a-1}q-\f{q-1}d=\f{(d^a-1)q+1}d>d^a,$$
also
$$2l<n+\f n{d^a}+\f{dn}q\ls n\l(1+\f1{d^2}+\f{d}{d+1}\r)\ls 2n.$$
So, we always have $1\ls k<l\ls n$ and hence we get a contradiction by the definition of $m$.
\medskip

{\it Case} 2. $q<d$.

If $q=1$, then $d^{h-1}<n\ls m=d^a\ls d^h$ and hence $m=d^h$ as desired.

Now suppose that $q>1$. As $q<d\ls 3$ we must have $q=2$ and $d=3$. Since $3^{h-1}<n\ls m=2\cdot 3^a\ls 3^h$,
we get $a=h-1$ and hence $3^a+1\ls n$. Observe that
$$(3^a+1)(3(3^a+1)-1)-1(3\cdot1-1)=3^a(3(3^a+2)-1)\eq0\ (\mo\ 2\cdot 3^a).$$
This contradicts that $m=2\cdot 3^a$.
\smallskip

Combining the above we have completed the proof of Theorem 1.2(ii). \qed

\proclaim{Lemma 3.2 {(\rm [RR, Theorem 1])}} Let $d\in\{1,\ldots,72\}$, and $r\in\Z$ with $(r,d)=1$.
For $x\gs10^{10}$ and $\ve=0.023269$, we have
$$(1-\ve)\f x{\varphi(d)}\ls \theta(x;r,d)\ls (1+\ve)\f x{\varphi(d)},$$
where $\varphi$ is Euler's totient function and
$\theta(x;r,d):=\sum\Sb p\ls x,\ p\eq r\, (\mo\ d)\endSb\log p$ with $p$ prime.
\endproclaim

\medskip
\noindent{\it Proof of Theorem} {\rm 1.2(iii)}. Let $n>3$ be an integer and take the least $m\in\Z^+$ with
$|\{18k(3k-1)\ \mo\ m:\ k=1,\ldots,n\}|=n$. We want to prove that $m$ is the least prime $p>3n$ with $p\eq1\ (\mo\ 3)$.
For $4\ls n\ls 36$ one can verify the desired result directly.

Below we assume $n>36$. Let $\ve=0.023269$. If $3n\gs 10^{10}$, then
$3.433(1-\ve)>3(1+\ve)$ and hence $\theta(3.433n;1,3)>\theta(3n;1,3)$ by Lemma 3.2,
therefore $(3n,3.433n]$ contains a prime $p\eq1\ (\mo\ 3)$.
For $n=37,\ldots,\lfloor10^{10}/3\rfloor$ one can easily verify (using a computer) that
the interval $(3n,3.433n]$ contains at least a prime congruent to 1 modulo 3.
(Note also that in 1932 R. Breusch [Br] refined the Bertrand Postulate confirmed by Chebyshev by showing that
for any $x\gs7$ the interval $(x,2x)$ contains a prime congruent to 1 modulo 3.)

If $p$ is a prime in $(3n,3.433n]$ with $p\eq1\ (\mo\ 3)$, then
for $1\ls k<l\ls n$ we have
$$18l(3l-1)-18k(3k-1)=18(l-k)(3(l+k)-1)\not\eq0\ (\mo\ p)$$
since $1\ls l-k<n<p$ and $p\not=3(l+k)-1<6n-1<2p$. Therefore $n\ls m\ls 3.433n$.

Assume that $m_0=m/(18,m)<3n$. As $m\gs n>36$ we have $m_0>2$. If $m_0\eq1\ (\mo\ 3)$, then
for $k=(m_0-1)/3$ and $l=(m_0+2)/3\ls n$ we have
$l(3l-1)\eq k(3k-1)\ (\mo\ m_0)$ and hence $18l(3l-1)\eq 18k(3k-1)\ (\mo\ m)$
which leads to a contradiction. As $4(3\cdot4-1)\eq 3(3\cdot3-1)\ (\mo\ 5)$, we cannot have $m_0=5$
since $k(3k-1)\ (k=1,\ldots,n)$ are pairwise distinct modulo $m_0$. If $m_0>5$ and $m_0\eq2\ (\mo\ 3)$, then
for $k=1<l=(m_0-2)/3\ls n$, we have $l(3l-1)\eq k(3k-1)\ (\mo\ m_0)$ which leads to a contradiction.
Therefore $3\mid m_0$. Write $m_0=3^aq$ with $a,q\in\Z^+$ and $3\nmid q$. If $q>1$, then we may argue as in cases 1 and 2 in the proof of
Theorem 1.2(ii) with $d=3$ to get a contradiction. So $m_0=3^a$,
and hence $m$ or $m/2$ is a power of 3. Suppose $3^{h-1}<n\ls 3^h$ with $h\in\Z^+$. Then
$m\in\{3^h,3^{h+1},2\cdot3^h,2\cdot3^{h-1}\}$ since $n\ls m\ls 3.433n$. For $k=1<l=3^{h-1}+1\ls n$ we clearly have
$m\mid 18(l-k)$ and hence $18l(3l-1)\eq 18k(3k-1)\ (\mo\ m)$ which leads to a contradiction.

By the above, we must have $m_0\gs 3n$. As $m/2<3n$ we must have $(18,m)=1$ and $m\gs 3n$.
If $p\in [3n,3.433n]$ is a prime with $p\eq2\ (\mo\ 3)$, then
for $k=(p-5)/6$ and $l=(p+7)/6$ we have $1\ls k<l\ls n$ and $18l(3l-1)\eq 18k(3k-1)\ (\mo\ p)$.

Now it remains to show that $m$ cannot be a composite number in $[3n,3.433n]$.
Suppose that $m=cd$ with $c,d\in\{2,3,\ldots\}$. As $(m,18)=1$, we have $(c,6)=(d,6)=1$.
Take $k\in[1,d]$ such that $k\eq ((1+2d(\f d3))/3-c)/2\ (\mo\ d)$, and set $l=k+c.$  Note that
$l(3l-1)-k(3k-1)=(l-k)(3(l+k)-1)\eq0\ (\mo\ m)$. Clearly
$$l=k+c\ls c+d=\f md+\f mc\ls 3.433n\l(\f1c+\f1d\r)\ls n$$
since $m=cd\gs n>36$ and $1/3.433\gs \max\{1/5+1/11,1/7+1/7\}$.  So we get a contradiction. \qed

\heading{4. Proofs of Theorems 1.3-1.5}\endheading
\proclaim{Lemma 4.1} Let $d\in\{4,6,12\}$ and $n\in\Z^+$. Then
$[2n-1,2.4n]$ contains at least a prime $p\eq -1\ (\mo\ d)$ except for $n\in E(d)$, where
$$E(4)=\{1,7,17\},\ \ \ E(6)=\{1,2,4,7,16,17\}$$
and
$$E(12)=\{1,2,3,4,7,8,9,13,14,15,16,17,18,19,43,44,67,68,69\}.$$
\endproclaim
\Proof. Note that $\ve:=0.023269<1/11$. If $n\gs 10^{10}/2$, then by Lemma 3.2 we have
$$\theta(2.4n;-1,d)\gs(1-\ve)\f{2.4n}{\varphi(d)}>(1+\ve)\f{2n}{\varphi(d)}\gs\theta(2n;-1,d),$$
and hence $(2n,2.4n]$ contains at least a prime $p\eq -1\ (\mo\ d)$.
It can be easily verified that for $n<10^{10}/2$ the interval $[2n-1,2.4n]$ contains a prime $p\eq-1\ (\mo\ d)$
except for $n\in E(d)$. We are done. \qed

\proclaim{Lemma 4.2} Suppose that $p>3$ is a prime in $[2n-1,2.4n]$ where $n>2$ is an integer. For $d\in\{4,6,12\}$,
those $(2k-1)^d$ with $1\ls k\ls n$ are pairwise distinct modulo $p$ if and only if $p\eq-1\ (\mo\ d)$.
\endproclaim
\Proof. For $1\ls k<l\ls n$, we clearly have
$$(2l-1)^4-(2k-1)^4=((2l-1)^2-(2k-1)^2)((2l-1)^2+(2k-1)^2),$$
$$\align&(2l-1)^6-(2k-1)^6=((2l-1)^3-(2k-1)^3)((2l-1)^3+(2k-1)^3)
\\=&((2l-1)^2-(2k-1)^2)((2l-1)^2+(2k-1)(2l-1)+(2k-1)^2)
\\&\times((2l-1)^2-(2k-1)(2l-1)+(2k-1)^2)
\endalign$$
and
$$\align&(2l-1)^6+(2k-1)^6
\\=&((2l-1)^2+(2k-1)^2)((2l-1)^4-(2k-1)^2(2l-1)^2+(2k-1)^4).
\endalign$$
Note that
$$(2l-1)^2-(2k-1)^2=4(l-k)(l+k-1)\not\eq0\pmod p$$
since $0<l-k<l+k-1<2n-1\ls p$. If $(2l-1)^2+(2k-1)^2\eq0\ (\mo\ p)$, then $-1$ is a quadratic residue mod $p$ and hence $p\eq1\ (\mo\ 4)$.
For $\delta\in\{\pm1\}$, if
$$4((2l-1)^2+\delta(2l-1)(2k-1)+(2k-1)^2)=(2(2l-1)+\delta(2k-1))^2+3(2k-1)^2$$
is divisible by $p$, then $-3$ is a quadratic residue mod $p$ and hence $p\eq1\pmod 6$.
Similarly, if $(2l-1)^4-(2k-1)^2(2l-1)^2+(2k-1)^4\eq0\ (\mo\ p)$ then $p\eq1\pmod 6$.

By the above, for any $d\in\{4,6,12\}$, if $p\eq-1\ (\mo\ d)$ then those $(2k-1)^d$ with $k=1,\ldots,n$ are pairwise distinct modulo $p$.

Now we handle the case $p\eq1\ (\mo\ 4)$. It is well known that $p=x^2+y^2$ for some integers $x>y>0$
and hence $2p=(x+y)^2+(x-y)^2$ with $x\pm y$ odd. Take $k=(x-y+1)/2$ and $l=(x+y+1)/2$. Clearly $2l-1=x+y\ls\sqrt{2p}\ls\sqrt{4.8n}<2n$
and hence $1\ls k<l\ls n$. As $(2l-1)^2\eq-(2k-1)^2\pmod p$, we have $(2l-1)^4\eq(2k-1)^4\ (\mo\ p)$ and $(2l-1)^{12}\eq(2k-1)^{12}\ (\mo\ p)$.

Now we assume $p\eq1\ (\mo\ 3)$. It is known that $p=u^2+3v^2$ for some $u,v\in\Z^+$ with $u\not\eq v\ (\mo\ 2)$.
Write $u+v=2l-1$ and $|u-v|=\delta(v-u)=2k-1$. Clearly $k,l\in\Z^+$ and $k<l$. Since $4p=(u-3v)^2+3(u+v)^2$, we have
$$u+v\ls\sqrt{\f{4p}3}\ls 2\sqrt{\f{2.4n}3}<2n$$
and hence $l\ls n$. Observe that
$$\align &(2l-1)^2+\delta(2l-1)(2k-1)+(2k-1)^2
\\=&(u+v)^2+(u+v)(v-u)+(u-v)^2=u^2+3v^2\eq0\pmod p.
\endalign$$
So we have $(2l-1)^6\eq (2k-1)^6\ (\mo\ p)$ and $(2l-1)^{12}\eq(2k-1)^{12}\ (\mo\ p)$.

Combining the above we have finished the proof of Lemma 4.2. \qed

\medskip
\noindent{\it Proof of Theorem} 1.3. Fix $d\in\{4,6,12\}$ and $n\in\{3,4,\ldots\}$.
We want to prove that $\lambda_d(n)$ (the least integer $m>1$ with $|\{(2k-1)^d\ \mo\ m:\ k=1,\ldots,n\}|=n$)
is just the least prime $p\gs2n-1$ with $p\eq-1\ (\mo\ d)$.

If $n\ls14$ or $n\in E(d)$, then we can easily verify the desired result. Below we simply assume $n\gs 15$ and $n\not\in E(d)$.

For $1\ls k<l\ls n$, clearly
$(2l-1)^d-(2k-1)^d$ is a multiple of $(2l-1)^2-(2k-1)^2=4l(l-1)-4k(k-1)$. If those $(2k-1)^d$ with $1\ls k\ls n$ are pairwise distinct
modulo an integer $m>1$, then so are those $k(k-1)$ $(k=1,\ldots,n)$ and hence $m\gs 2n-1$ by Theorem 2.1(i). Therefore $\lambda_d(n)\gs 2n-1$.

By Lemma 4.1, $[2n-1,2.4n]$ contains a prime $p\eq-1\ (\mo\ d)$ and hence $\lambda_d(n)\ls p\ls 2.4n$ by Lemma 4.2.
As those $2k(k-1)\ (k=1,\ldots,n)$ are pairwise distinct mod $\lambda_d(n)$, by Theorem 2.1(ii) and Theorem 2.2(ii),
$\lambda_d(n)$ must be a prime. In view of Lemma 4.2, $\lambda_d(n)$ is the least prime $p\in[2n-1,2.4n]$ with $p\eq-1\ (\mo\ d)$.

So far we have completed the proof of Theorem 1.3. \qed

\proclaim{Lemma 4.3} For any odd prime $q$ and positive integer $n$, the interval $[2n-1,2.4n]$ contains at least a prime $p\not\eq1\ (\mo\ q)$
unless $n\ls 17$ and $q<2.4n$.
\endproclaim
\Proof. By the proof of Lemma 3.1,  $[2n-1,2.4n]$ contains a prime $p$. If $p\eq1\ (\mo\ q)$ then $q\ls p-1<2.4n$.

Clearly $[2\cdot1-1,2.4]$ contains the prime $2\not\eq1\ (\mo\ q)$.
When $n>1$, the interval $[2n-1,2.4n]$ contains an odd prime $p$. If $q\gs 1.2n$ then $1+2q>2.4n$ and hence $p\not\eq1\ (\mo\ q)$.
Below we assume $q<1.2n$.

We first handle the case $q\ls 53$. As in Lemma 4.1 we can employ [RR, Theorem 1.1] to deduce that $(2n,2.4n]$ contains a prime $p\eq-1\ (\mo\ q)$
for $n\gs 10^{10}/2$. For $n\in[18,10^{10}/2]$ we can easily check that $[2n-1,2.4n]$ indeed contains a prime $p\not\eq1\ (\mo\ q)$.

Now assume that $q\gs 59$. Set $x:=2.4n$. Then $q<x/2$.
By the Brun-Titchmarsh theorem (cf. [MV] or [CP, p.\,43]) in analytic number theory, we have
$$\pi(x;1,q):=|\{p\ls x:\ p\ \t{is a prime and}\ p\eq1\ (\mo\ q)\}|\ls\f{2x}{\varphi(q)\log(x/q)}.$$
Thus, if $q\ls\sqrt x$ then
$$\pi(x;1,q)\ls\f{2x}{(q-1)\log\sqrt x}\ls\f{4x}{58\log x}=\f2{29}\times\f{x}{\log x};$$
if $\sqrt x<q\ls x/2$ then
$$\pi(x;1,q)\ls\f{2x}{(\sqrt x-1)\log2}.$$
Note that $(\sqrt x-1)\log 2>29\log x$ when $n\gs 114895$.

Assume $n>148000$. By the above,
$$\pi(x;1,q)\ls \f 2{29}\times\f{x}{\log x}.\tag4.1$$
Since $x=2.4n>599$, by [D, Section 4] we have
$$\pi(x):=\pi(x;1,1)\gs \f x{\log x}\l(1+\f{0.992}{\log x}\r)>\f x{\log x}$$
and
$$\align\pi(2n)\ls &\f {2n}{\log (2n)}\l(1+\f{1.2762}{\log (2n)}\r)
\\\ls&\f{2n}{\log(2n)}\l(1+\f{1.2762}{\log(2\times148001)}\r)<\f{2.202602n}{\log(2n)}.
\endalign$$
Thus
$$\pi(2.4n)-\pi(2n)
>\f{2.4n}{\log(2.4n)}-\f{2.202602n}{\log(2n)}.\tag4.2$$
Since
$$\align \l(\f{27}{29}\times 2.4-2.202602\r)\log n&\gs\l(\f{27}{29}\times 2.4-2.202602\r)\log 148001
\\&>0.3795>2.202602\log2.4-\f{27}{29}\times 2.4\log 2,
\endalign$$
we have the inequality
$$\l(1-\f 2{29}\r)\f{2.4}{\log n+\log 2.4}>\f{2.202602}{\log n+\log 2}.\tag4.3$$
Combining (4.1)--(4.3) we obtain $\pi(2.4n)-\pi(2n)>\pi(2.4n;1,q)$. So $[2n-1,2.4n]$ contains a prime $p\not\eq1\ (\mo\ q)$.

When $18\ls n\ls 148000$ and $59\ls q<1.2n$, we can easily verify the desired result using a computer.

So far we have proved Lemma 4.3.
\qed

\medskip
\noindent{\it Proof of Theorem} {\rm 1.4}. Fix an odd prime $q$ and let $D_q(n)$ denote the smallest integer $m>1$ such that
those $k^q(k-1)^q$ $(k=1,\ldots,n)$ are pairwise distinct modulo $m$. We want to prove that $D_q(n)$ is just the least prime $p\gs2n-1$
with $p\not\eq1\ (\mo\ q)$. This is trivial for $n=1$, so we just let $n>1$.

As those $k(k-1)\ \mo\ D_q(n)$ with $1\ls k\ls n$
are pairwise distinct, we have $2<2n-1\ls T(n)\ls D_q(n)$ by Theorem 1.1(ii).

If $n\ls 17$ and $q<2.4n$, then we can easily verify the desired result directly.
Below we let $n\gs 18$ or $q\gs2.4n$. By Lemma 4.3, the interval $[2n-1,2.4n]$ contains a prime $p\not\eq1\ (\mo\ q)$.

Let $p$ be any prime in $[2n-1,2.4n]$. If $l^q(l-1)^q\eq k^q(k-1)^q\ (\mo\ p)$ for some $1\ls k<l\ls n\ls (p+1)/2$,
then $p\nmid k(k-1)$,
$$\l(\f{l(l-1)}{k(k-1)}\r)^q\eq 1\ (\mo\ p)\ \ \t{and}\ \l(\f{l(l-1)}{k(k-1)}\r)^{(q,p-1)}\eq 1\ (\mo\ p);\tag4.4$$
as $l(l-1)\not\eq k(k-1)\ (\mo\ p)$ by Theorem 2.2(i), (4.4) implies that $(q,p-1)>1$ and hence $p\eq1\ (\mo\ q)$.
Conversely, if $p\eq1\ (\mo\ q)$, then $q<p\ls 2.4n$ and $n\gs18$, hence those $k^q(k-1)^q$ with $1\ls k\ls n$ cannot be pairwise distinct
modulo $p$ since we only have $(p-1)/q\ls (p-1)/3\ls(2.4n-1)/3<n-1$  $q$-th power residue modulo $p$.

In view of the above, $D_q(n)$ does not exceed the least prime $p\in[2n-1,2.4n]$ with $p\not\eq1\ (\mo\ q)$.
If $D_q(n)=2^aw$ with $a\gs 3$ and $2\nmid w$, then
$$(2^{a-2}w(2^{a-2}w-1))^q\eq(1(1-1))^q\ \ (\mo\ 2^aw)$$
and also $1<2^{a-2}w=D_q(n)/4\ls0.6n<n$. So $8\nmid D_q(n)$. If $D_q(n)=2^aw$ with $a\in\{1,2\}$ and $2\nmid w$, then
those $k^q(k-1)^q$ $(k=1,\ldots,n)$ are pairwise distinct modulo $w<D_q(n)$ since $8\mid k^q(k-1)^q$ for all $k=1,\ldots,n$.
Thus $D_q(n)$ cannot be even.
If $n\gs 15$, then $D_q(n)$ must be a prime
by Theorem 2.1(ii), and hence it is just the least prime $p\gs 2n-1$ with $p\not\eq1\ (\mo\ q)$.

Now we handle the remaining case $2\ls n\ls 14$ and $q\gs2.4n$. Note that any prime in $[2n-1,2.4n]$
is not congruent to 1 modulo $q$. For each $n=2,\,3,\,4,\,6,\,7,\,9,\,10,\,12$, clearly $2n-1$ is prime and hence $D_q(n)$ is the least prime
in $[2n-1,2.4n]$. As $2+9/3=5$, we have $3^2\nmid D_q(5)$ by the proof of Lemma 2.2, hence $D_q(5)$ is the least prime $11$ after $2\cdot5-1=9$.
(Note that $D_q(5)\not=10$ since 10 is even.) Since $15/3+3=8$ and $21/3+3<11$, by the proof of Theorem 2.1(ii) we have
$D_q(8)\not=3\cdot 5$ and $D_q(11)\not=3\cdot7$, hence $D_q(8)=17$ and $D_q(11)=23$ as desired. For $n=13,14$,
as $2+D_q(n)/3\ls 2+0.8n\ls n$, by the proof
of Lemma 2.2 we have $p^2\nmid D_q(n)$ for any odd prime $p$, hence $D_q(n)\not=25,\,27$. Note also that $D_q(n)\not=26,\,28$. So
$D_q(13)=D_q(14)=29$ as desired.

The proof of Theorem 1.4 is now complete. \qed

\proclaim{Lemma 4.4} All those $s_n=\sum_{k=1}^n(-1)^{n-k}p_k\ (n=1,2,3,\ldots)$ are pairwise distinct, and also $s_n\ls p_n$ for all $n\in\Z^+$.
\endproclaim
\Proof. Obviously $s_1=p_1=2$. For $n=2,3,4,\ldots$,  we clearly have $s_n+s_{n-1}=p_n$ and hence $s_n<p_n$ since $s_{n-1}>0$.

Now we show that $s_n\not=s_k$ for any $1\ls k<n$ (see also [SC] for this simple observation). If $n-k$ is even, then
$$s_n-s_k=(p_n-p_{n-1})+\cdots+(p_{k+2}-p_{k+1})>0.$$
When $n-k$ is odd, we have
$$s_n-s_k=\sum_{l=k+1}^n(-1)^{n-l}p_l-2\sum_{j=1}^k(-1)^{k-j}p_j\eq n-k\not\eq0\pmod2.$$

The proof of Lemma 4.4 is now complete. \qed

\medskip
\noindent{\it Proof of Theorem 1.5}. Let $k,l\in\{1,\ldots,n\}$ with $k\not=l$. We want to show that
$$2s_l^2-2s_k^2=2(s_l+s_k)(s_l-s_k)\not\eq0\pmod{p_{n+1}}.$$

By Lemma 4.4, $s_k\not=s_l$ and $|s_k-s_l|\ls\max\{s_k,s_l\}\ls \max\{p_k,p_l\}\ls p_n<p_{n+1}$,
therefore $s_k\not\eq s_l\ (\mo\ p_{n+1})$.

As $s_k+s_l\ls p_k+p_l\ls 2p_n<2p_{n+1}$, it remains to prove that $s_k+s_l\not=p_{n+1}$.
Without loss of generality we assume that $k<l$. If $l-k$ is even, then
$$s_l+s_k=\sum_{j=k+1}^l(-1)^{l-j}p_j+2s_k\eq l-k\eq0\pmod{2}$$
and hence $s_k+s_l\not=p_{n+1}$. If $l-k$ is odd, then
$$s_l+s_k=\sum_{j=k+1}^l(-1)^{l-j}p_j=p_l-\sum_{0<j\ls(l-k-1)/2}(p_{l-2j+1}-p_{l-2j})\ls p_l\ls p_n<p_{n+1}.$$
So we do have $s_k+s_l\not=p_{n+1}$ as desired.

In view of the above we have completed the proof of Theorem 1.5. \qed

\heading{5. More conjectures}\endheading

Motivated by Conjecture 1.1, here we pose more conjectures for further research.

\proclaim{Conjecture 5.1} {\rm (i)} For the functions $s(n)$ and $t(n)$ in Conjecture {\rm 1.1}, we have
$s(n)<n^2$ and $t(n)\ls n^2/2$ for all $n=2,3,4,\ldots$.

{\rm (ii)} The number of primes not exceeding $x$ in the set
$S=\{s(1),s(2),s(3),\ldots\}$ is $o(\sqrt x)$ and even $O(\sqrt x/\log^3 x)$ as $x\to+\infty$.

{\rm (iii)} If we replace $k!$ in Conjecture {\rm 1.1(ii)} by $(k+1)!$ or $(2k)!$, then the modified $t(n)$ is always a prime.
\endproclaim
\Remark\ 5.1. It seems that if we replace $\bi{2k}k$ in the definition of $s(n)$ by $2^{k!}$ or $2^k!$ or $2^{2^k}$
then the modified $s(n)$ also takes only prime values.
\medskip

\proclaim{Conjecture 5.2} Let $n$ be a positive integer.

{\rm (i)} The least integer $m>1$ such that $|\{(k^2-k)!\ \mo\ m:\ k=1,\ldots,n\}|=n$ is a prime
in the interval $((n-1)(n-2),n(n-1))$ for every $n=3,4,\ldots$.

{\rm (ii)} The least integer $m>1$ such that $n!\not\eq k!\ (\mo\ m)$
for all $0<k<n$ is a prime not exceeding $2n$ except for $n=4,6$.
\endproclaim
\Remark\ 5.2. For any positive integer $n$, the interval $[n,2n]$ contains at least a prime
by the Bertrand Postulate proved by Chebyshev, but Legendre's conjecture that $(n^2,(n+1)^2)$ contains a prime remains unsolved.
\medskip

\proclaim{Conjecture 5.3} Let $a\in\Z$ with $|a|>1$. For $n\in\Z^+$ define $f_a(n)$ as the least integer $m>1$ such that
those $a^k\ (k=1,\ldots,n)$ are pairwise distinct modulo $m$. Then there is a positive integer $n_0(a)$ such that
for any integer $n\gs n_0(a)$, the number $f_a(n)$ is the least prime $p>n$ having $a$ as a primitive root modulo $p$
if $a$ is not a square, and  $f_a(n)$ is the least prime $p>2n$  such that
$a,a^2,\ldots,a^{(p-1)/2}$ are pairwise distinct modulo $p$ if $a$ is a square.
In particular, we may take
$n_0(-2)=3,\ n_0(-3)=n_0(5)=1$, and $n_0(9)=n_0(25)=2$.
\endproclaim

Let $A$ and $B$ be integers. The Lucas sequence $u_n=u_n(A,B)\ (n\in\N=\{0,1,2,\ldots\})$
and its companion sequence $v_n=v_n(A,B)\ (n\in\N)$ are defined as follows:
$$u_0=0,\ u_1=1,\ \t{and}\ u_{n+1}=Au_n-Bu_{n-1}\ (n=1,2,3,\ldots);$$
and
$$v_0=2,\ v_1=A,\ \t{and}\ v_{n+1}=Av_n-Bv_{n-1}\ (n=1,2,3,\ldots).$$
It is well known that
$$(\al-\beta)u_n=\al^n-\beta^n\ \ \t{and}\ \ v_n=\al^n+\beta^n\quad\t{for all}\ n\in\N,$$
where $\al=(A+\sqrt{\Delta})/2$ and $\beta=(A-\sqrt{\Delta})/2$ are the two roots of the equation $x^2-Ax+B=0$ with $\Delta=A^2-4B$.
It is also known that if $p$ is an odd prime not dividing $B$ then  $p\mid u_{p-(\f{\Delta}p)}$ (see, e.g., [S06]),
where $(-)$ is the Legendre symbol.
Note that
$$u_{2n}=u_nv_n=Au_n(A^2-2B,B^2)\ \ \t{and}\ \ v_{2n}=v_n(A^2-2B,B^2)$$ for all $n\in\N$.
Those $F_n=u_n(1,-1)$ and $L_n=v_n(1,-1)$ are Fibonacci numbers and Lucas numbers respectively, and also
$F_{2n}=u_n(3,1)$ and $L_{2n}=v_n(3,1)$.

Clearly an integer $a$ is a primitive root modulo an odd prime $p$ if and only if those $v_k(a+1,a)=a^k+1\ (k=1,...,p-1)$
 are pairwise distinct modulo $p$. Motivated by the Artin conjecture, we raise the following new conjecture.

\proclaim{Conjecture 5.4} Let $A$ be an integer with $|A|>2$.

  {\rm  (i)} If $2+A$ is not a square, then there are infinitely many odd primes $p\nmid A^2-4$ such that those $v_k(A,1)\ \mo\ p$
   for $k=1,...,(p-(\f{A^2-4}p))/2$ are pairwise distinct.

  {\rm (ii)} If $2-A$ is not a square, then there are infinitely many odd primes $p\nmid A^2-4$ such that those
  $u_k(A,1)\ \mo\ p$ for $k=1,...,(p-(\f{A^2-4}p))/2$ are pairwise distinct.
\endproclaim

Inspired by Conjecture 5.3, we pose the following challenging conjecture which implies part (i) of Conjecture 5.4.

\proclaim{Conjecture 5.5} Let $A$ be an integer with $|A|>2$. For $n\in\Z^+$ define $t_A(n)$ as the smallest integer $m>1$
such that those $v_k(A,1)\ \mo\ m$ for $k=1,\ldots,n$ are pairwise distinct.
Then  $t_A(n)$ is prime for any sufficiently large integer $n$
$(n>2|A|$ might suffice$)$. When $A+2$ is not a square, there is a positive integer $N_0(A)$ such that
for any integer $n\gs N_0(A)$, the number $t_A(n)$ is the smallest odd prime $p\nmid A^2-4$ such that
$p-(\f{A^2-4}p)\gs2n$ and those $v_k(A,1)\ \mo\ p$ $(k=1,\ldots,(p-(\f{A^2-4}p))/2)$ are pairwise distinct.
In particular, we may take $N_0(3)=6$, $N_0(-3)=7$, and $N_0(\pm4)=N_0(\pm10)=3$.
 \endproclaim

\Remark\ 5.3. Note that $v_k(3,1)=L_{2k}$ and $v_k(-3,1)=(-1)^kL_{2k}$ for any $k\in\Z^+$.
Also, [S02] contains the congruence
$$T_{(p-(\f 3p))/2}\eq 2\l(\f 6p\r)\pmod{p^2}\ \quad\t{for any prime}\ p>3,$$
where $T_n:=v_n(4,1)$.

\medskip

Recall that $S_n$ denotes the sum of the first $n$ primes. Our following conjecture is a refinement of the Artin conjecture.

\proclaim{Conjecture 5.6}
If $a\in\Z$ is neither $-1$ nor a square, then there is a positive integer $n_0$ such that for any integer $n\gs n_0$
the least integer $m>1$ such that $|\{a^{S_k}\ \mo\ m:\ k=1,\ldots,n\}|=n$
is a prime $p$ having $a$ as a primitive root modulo $p$.
In particular, we may take $n_0=1$ for $a=-3$.
\endproclaim

Recall that the Euler numbers $E_0,E_1,E_2,\ldots$ are integers defined by
$$E_0=1,\ \ \t{and}\ \ \sum^n\Sb k=0\\2\mid k\endSb \bi nk E_{n-k}=0\ \ \ \t{for}\ n=1,2,3,\ldots.$$
It is well known that $E_{2n+1}=0$ for all $n\in\N$ and
$$\sec x=\sum_{n=0}^\infty(-1)^n E_{2n}\f{x^{2n}}{(2n)!}\ \ \l(|x|<\f{\pi}2\r).$$

\proclaim{Conjecture 5.7} {\rm (i)} For $n\in\Z^+$ let $e(n)$ be the least integer $m>1$ such that $E_{2k}\ (k=1,\ldots,n)$
are pairwise distinct modulo $m$. Then we have $e(n)=2^{\lceil\log_2n\rceil+1}$ with the only exceptions as follows:
$$\align &e(3)=7,\ e(5)=e(6)=13,\ e(9)=e(10)=25,\ e(17)=47,
\\&e(18)=e(19)=e(20)=e(21)=7^2,\ \ e(65)=\cdots=e(78)=13^2,
\\&e(1025)=e(1026)=e(1027)=e(1028)=e(1029)=e(1030)=5^5.
\endalign$$

{\rm (ii)} For $n\in\Z^+$ let $e^*(n)$ be the least integer $m>1$ such that $2E_{2n}\not\eq 2E_{2k}\ (\mo\ m)$ for all $0<k<n$.
Then $e^*(n)$ is a prime in the interval $[2n,3n]$ with the only exceptions as follows:
$$e^*(4)=13,\ e^*(7)=23,\ e^*(10)=5^2,\ e^*(55)=11^2.$$
\endproclaim
\Remark\ 5.4. With the help of the Stern congruence for Euler numbers (see, e.g., S. S. Wagstaff [W] and the author [S05]),
we can easily show that $\log_2 e(n)\ls \lceil\log_2n\rceil+1$. Also, it is known (cf. [B]) that for any $n\in\Z^+$ the interval
$[2n,3n]$ contains at least a prime.
\medskip

\Ack. The author would like to thank Prof. N. Koblitz, C. Pomerance, P. Moree,
and Dr. O. Gerard and H. Pan, and the referee for their helpful comments.


\widestnumber\key{ABM}

 \Refs

\ref\key ABM\by L. K. Arnold, S. J. Benkoski and B. J. McCabe\paper The discriminator (a simple application of Bertrand's postulate)
\jour Amer. Math. Monthly\vol 92\yr 1985\pages 275--277\endref

\ref\key B\by M. El Bachraoui\paper Primes in the interval $[2n,3n]$\jour Int. J. Contemp. Math. Sci.\vol 1\yr 2006\pages 617--621\endref

\ref\key BSW\by P. S. Bremser, P. D. Schumer and L. C. Washington
\paper A note on the incongruence of consecutive integers to a fixed power
\jour J. Number Theory\vol 35\yr 1990\pages 105--108\endref

\ref\key Br\by R. Breusch\paper Zur Verallgemeinerung des Bertrandschen Postulates, dass zwischen $x$ und $2x$ stets Primzahlen liegen
\jour Math. Z.\vol 34\yr 1932\pages 505--526\endref

\ref\key CP\by R. Crandall and C. Pomerance\book Prime Numbers: A Computational Perspective\publ 2nd Edition, Springer, New York, 2005\endref

\ref\key D\by P. Dusart\paper The $k$th prime is greater than $k(\log k+\log\log k-1)$ for $k\gs2$
\jour Math. Comp. \vol 68\yr 1999\pages 411--415\endref

\ref\key H\by W. B. Hart\paper Re: A new conjecture on primes\jour a message to Number Theory List, April 14, 2012.
{\tt http://listserv.nodak.edu/cgi-bin/wa.exe?A2=NMBRTHRY;57b2e5f8.1204}\endref

\ref\key IR\by K. Ireland and M. Rosen\book A Classical Introduction to Modern Number Theory, 2nd Edition
\publ Springer, New York, 1990\endref

\ref\key Mi\by W. H. Mills\paper A prime-representing function\jour Bull. Amer. Math. Soc.\vol 53\yr 1947\pages 604\endref

\ref\key MV\by H. Montgomery and R. Vaughan\paper The large sieve\jour Mathematica\vol 20\yr 1973\pages 119--134\endref

\ref\key MM\by P. Moree and G. L. Mullen\paper Dickson polynomial discriminators\jour J. Number Theory \vol 59\yr 1996\pages 88--105\endref

\ref\key Mu\by R. Murty\paper Artin's conjecture for primitive roots\jour Math. Intelligencer\vol 10\yr 1988\pages 59--67\endref

\ref\key N\by J. Nagura\paper On the interval containing at least one prime number\jour Proc. Japan Acad. Ser. A\vol 28\yr 1952\pages 177--181\endref

\ref\key RR\by O. Ramar\'e and R. Rumely\paper Primes in arithmetic progressions\jour Math. Comp.
\vol 65\yr 1996\pages 397--425\endref

\ref\key SC\by N.J.A. Sloane and J. H. Conway\paper {\rm Sequence A008347 in OEIS
(On-Line Encyclopedia of Integer Sequences)} \jour {\tt http://oeis.org/A008347}\endref

\ref\key S02\by Z. W. Sun\paper {\rm On the sum $\sum_{k\eq r\,(\mo\ m)}\bi nk$
and related congruences}\jour Israel J. Math.
\vol 128\yr 2002\pages135--156\endref

\ref\key S05\by Z. W. Sun\paper On Euler numbers modulo powers of two\jour J. Number Theory\vol 115\yr 2005\pages 371--380\endref

\ref\key S06\by Z. W. Sun\paper Binomial coefficients and quadratic fields
\jour Proc. Amer. Math. Soc.\vol 134\yr 2006\pages 2213--2222\endref

\ref\key S12a\by Z. W. Sun\paper A function taking only prime values\jour a message to Number Theory List on Feb. 21, 2012.
{\tt http://listserv.nodak.edu/cgi-bin/wa.exe?A2=ind1202\&L=NMBRTHRY}
{\tt \&T=0\&P=4035}\endref

\ref\key S12b\by Z. W. Sun\paper {\rm Sequence A208494 in OEIS
(On-Line Encyclopedia of Integer Sequences), posted in Feb. 2012}
\jour {\tt http://oeis.org/A208494}\endref

\ref\key S12c\by Z. W. Sun\paper {\rm Sequences A210144 and A210186 in OEIS
 (On-Line Encyclopedia of Integer Sequences)}, {\rm posted in March 2012}
\jour {\tt http://oeis.org}\endref

\ref\key S13\by Z. W. Sun\paper {\rm Sequences A222579 and A222580 in OEIS
 (On-Line Encyclopedia of Integer Sequences)}, {\rm posted in Feb. 2013}
\jour {\tt http://oeis.org}\endref


\ref\key W\by S. S. Wagstaff, Jr.\paper Prime divisors of the Bernoulli and Euler numbers
\jour in: Number Theory for the Millennium, III (Urbana, IL, 2000), 357--374,
A K Peters, Natick, MA, 2002\endref

\ref\key Z\by M. Zieve\paper A note on the discriminator\jour J. Number Theory\vol 73\yr 1998\pages 122--138\endref


\endRefs
\enddocument